\documentstyle[11pt]{article}
\begin{document}
\title{\textbf{Reflected Solutions of Backward Doubly Stochastic \\Differential Equations\footnote{Supported by National Natural Science Foundation of China Grant 10771122,
Natural Science Foundation of Shandong Province of China Grant
Y2006A08 and National Basic Research Program of China (973 Program,
No. 2007CB814900)}}\\}
\author{Yufeng SHI \thanks{
Corresponding author. E-mail: yfshi@sdu.edu.cn},\quad Yanling GU\\
School of Mathematics \\
Shandong University\\
Jinan, 250100, P. R. China}
\date{}
\maketitle
\begin{abstract}
We study reflected solutions of one-dimensional backward doubly
stochastic differential equations (BDSDEs in short). The
``reflected'' keeps the solution above a given stochastic process.
We get the uniqueness and existence by penalization. For the
existence of backward stochastic integral, our proof is different
from [KKPPQ] slightly. We also obtain a comparison theorem for
reflected BDSDEs. At last we gave a simulation for the reflected
solutions of BDSDEs.
\end{abstract}


\section{ Introduction}
\section{Preliminaries: the existence and uniqueness to BDSDEs}
\textbf{Notations.} The Euclidean norm of a vector
$x\in{\mathbf{R}}^k$ will be denoted by $|x|$, and for a
$d\times k$ matrix $A$, we define $\|A\|=\sqrt{TrAA^*}$.\\
Let $(\Omega,\cal{F},P)$ be a probability space, and $T>0$ be an
arbitrarily fixed constant throughout this paper. Let $\{W_t;0\leq
t\leq T\}$ and $\{B_t;0\leq t\leq T\}$ be two mutually independent
standard Brownian Motion processes, with values respectively in
$R^d$ and $R^l$, defined on $(\Omega,\cal{F},P)$. Let $\cal{N}$
denote the class of $P-$null sets of $\cal{F}$. For each $t\in
[0,T]$, we define
$${\cal{F}}_t : ={\cal{F}}_t^W\vee {\cal{F}}_{t,T}^B$$
where for any process $\{\eta_t\}$,
${\cal{F}}_{s,t}^{\eta}=\sigma\{\eta_r-\eta_s;s\leq r\leq t\}\vee
{\cal{N}}$, ${\cal{F}}_t^{\eta}={\cal{F}}_{0,t}^{\eta}$.

Note that the collection $\{{\cal{F}}_t;t\in [0,T]\}$ is neither
increasing nor decreasing, so it does not constitute a
filtration.\\
Let us introduce some spaces.

${\mathbf{L}}^2=\{\xi$ is an ${\cal{F}}_T$ -measurable random
variable s.t. $E(|\xi|^2)<\infty \}$.

${\mathbf{M}}^n=\{\{\phi_t,0\leq t\leq T\}$ is a jointly
measurable (classes of $dP\times dt$ a.e. equal) stochastic
processes s.t. $E\int_0^T|\varphi_t|^ndt<\infty$, and $\phi_t$ is
${\cal{F}}_t$ measurable for a.e. $t\in [0,T]\}$. ($n\in
{\mathbf{N}}$)

${\mathbf{S}}^2=\{\{\phi_t,0\leq t\leq T\}$ is a continuous
stochastic processes, s.t. $E(\sup_{0\leq t\leq
T}|\phi_t|^2)<\infty$, and
$\phi_t$ is ${\cal{F}}_t$ measurable for a.e. $t\in [0,T]\}$.\\
We are given four objects: the first one is a terminal value
$\xi$, s.t.

(i) $\xi \in {\mathbf{L}^2}$.\\
The second is two ``coefficients'' $f$ and $g$, which are maps
\begin{eqnarray*}
f:\Omega\times [0,T]\times R\times R^d\rightarrow R, \quad \quad
g:\Omega\times [0,T]\times R\times R^d\rightarrow R^l.
\end{eqnarray*}
 be jointly measurable and such that

(ii) $\forall (y,z)\in R\times R^d$, $f(\cdot,y,z)\in
{\mathbf{M}}^2$, $g(\cdot,y,z)\in {\mathbf{M}}^2$.

(iii) there exist two constants $K>0$ and $0<\alpha<1$, $\forall
y$, $y'\in R$, $z$, $z'\in R^d$, a.s.,a.e.
$$
\left\{
\begin{array}{lcl}
|f(t,y,z)-f(t,y',z')|^2 &\leq & C(|y-y'|^2+|z-z'|^2)\\
|g(t,y,z)-g(t,y',z')|^2 &\leq &C|y-y'|^2+\alpha
|z-z'|^2\\
\end{array}
\right. \eqno{\mbox{(H1)}}
$$
And the last one is an ``obstacle'' $\{S_t,0\leq t\leq T\}$, which
is a continuous progressively measurable real-valued process,
$S_t$ is ${\cal{F}}_t$ measurable, satisfying

(iv) $E\{\sup_{0\leq t\leq T}(S^+_t)\}<\infty$.\\
We always assume that $S_T\leq \xi$, a.s.\\
Now, let us introduce our reflected backward doubly stochastic
differential equations (BDSDE in short). The solution of our
reflected BDSDE is triple $(Y,Z,K)$ of ${\cal{F}}_t$ measurable
processes taking valued in $\bf{R}$, $\bf{R}^d$ and $\bf{R}_+$,
respectively, and satisfying

(v) $Z\in {\mathbf{M}}^2$;

(v') $Y\in {\mathbf{S}}^2$, and $K_T\in {\mathbf{L}}^2$;

(vi)
$Y_t=\xi+\int_t^Tf(s,Y_s,Z_s)ds+K_t-K_t+\int_t^Tg(s,Y_s,Z_s)dB_s-\int_t^TZ_s,dW_s$,
$0\leq t\leq T$;

(vii) $Y_t\geq S_t$, \quad\quad $0\leq t\leq T$;

(viii) $\{K_t\}$ is continuous and increasing, $K_0=0$ and
$\int_0^T(Y_t-S_t)dK_t=0$.

\textbf{Lemma 2.1.} Under the above conditions (i), (ii) and
(iii), the following backward doubly stochastic differential
equation (BDSDE in short) (1)
\begin{eqnarray}
Y_t=\xi+\int_t^Tf(s,Y_s,Z_s)ds+\int_t^Tg(s,Y_s,Z_s)dB_s-\int_t^TZ_sdW_s,\quad
0\leq t\leq T.
\end{eqnarray}
has a unique solution $(Y,Z)\in {\mathbf{S}}^2\times {\mathbf{M}}^2.$\\
This lemma was derived from Pardoux and Peng [PP3].

\textbf{Lemma 2.2.} let $(\xi,f,g)$ and $(\xi',f',g)$ be two
parameters of BDSDEs, each one satisfies all the assumptions (i),
(ii) and (iii) [with the exception that the Lipschitz condition
(iii) could be satisfied by either $f$ or $f'$ only], and suppose
in addition the following
$$\xi\leq \xi', a.s., \quad \quad f(t,y,z)\leq f'(t,y,z), a.s.a.e \quad \forall (y,z)\in R\times R^d .$$
Let $(Y,Z)$ be a solution of the BDSDE with parameter $(\xi,f,g)$
and $(Y',Z')$ a solution of the BDSDE with parameter
$(\xi',f',g)$. Then
$$Y_t\leq Y'_t, \quad a.e. \quad \quad \forall  0\leq t\leq T$$
This lemma was derived from Y. Gu and Y. Shi [GS].
\section{A priori estimate.}
In the following, $c$ will denote a constant whose value can vary
from line to line. \\
{\bf Proposition 3.1.} let $(Y,Z,K)$ be a solution of the
following reflected BDSDE(2)
\begin{eqnarray}
Y_t=\xi+\int_t^Tf(s,Y_s,Z_s)ds+\int_t^Tg(s,Y_s,Z_s)dB_s+K_t-K_t-\int_t^TZ_sdW_s,\quad
0\leq t\leq T.
\end{eqnarray}
Then there exists a constant $C$ such that
\begin{eqnarray}\nonumber
& & E(\sup_{0\leq t\leq T}Y_t^2+\int_0^T|Z_t|^2dt+K_T^2)\\
&\leq& CE(\xi^2+\int_0^ Tf(t,0,0)^2dt+\int_0^
Tg(t,0,0)^2dt+\sup_{0\leq t\leq T}(S_t^+)^2).
\end{eqnarray}
{\bf Proof.} Applying It$\hat{o}$'s to the process $Y_t$ and the
function $y\rightarrow y^2$ yields,
\begin{eqnarray*}
Y_t^2+\int_t^T|Z_s|^2ds&=&\xi^2+2\int_t^TY_sf(s,Y_s,Z_s)ds+\int_t^Tg^2(s,Y_s,Z_s)ds\\
& &
+2\int_t^TS_sdK_s+2\int_t^TY_sg(s,Y_s,Z_s)dB_s-2\int_t^TY_sZ_sdW_s
\end{eqnarray*}
where we have used the identity $\int_0^T(Y_t-S_t)dK_t=0$.\\
Using the Lipschitz property of $f$ and $g$, we have
\begin{eqnarray*}
EY_t^2+E\int_t^T|Z_t|^2 dt&=&
E\xi^2+2E\int_t^T|Y_s|\cdot[K(|Y_s|+|Z_s|)+f(s,0,0)]ds+2E\int_t^TS_sdK_s\\
& &
+E\int_t^Tc(\alpha')(|Y_s|^2+|g(s,0,0)|^2)ds+\alpha'E\int_t^T|Z_s|^2ds\\
&\leq&
E\xi^2+(2K+1+K^2\frac{2}{1-\alpha'}+c(\alpha'))E\int_t^T|Y_s|^2ds+E\int_t^Tf^2(s,0,0)ds\\
& &
+c(\alpha')E\int_t^Tg^2(s,0,0)ds+2\int_t^TS_sdK_s+\frac{1+\alpha'}{2}E\int_t^T|Z_s|^2ds
\end{eqnarray*}
where we have used (iii), for any $\alpha<\alpha'<1$, $\exists
c(\alpha')$, such that
\begin{eqnarray} |g(t,y,z)|^2\leq
c(\alpha')(|y|^2+|g(t,0,0)|^2)+\alpha'|z|^2, \quad a.s.\quad
\forall (y,z)\in R\times R^d
\end{eqnarray}
from Gronwall's lemma applied to $Y$, we get
\begin{eqnarray}
EY_t^2\leq
cE[\xi^2+\int_t^Tf(s,0,0)^2ds+\int_t^T|g(s,0,0)|^2ds+2\int_t^TS_sdK_s]
\end{eqnarray}
It follows that
\begin{eqnarray}
E\int_0^T|Z_s|^2ds\leq
cE[\xi^2+\int_0^Tf(s,0,0)^2ds+\int_0^T|g(s,0,0)|^2ds+2\int_0^TS_sdK_s]
\end{eqnarray}
We now give an estimate of $E[K_T^2]$, from Eq(2)
\begin{eqnarray*}
K_T=Y_0-\xi-\int_0^Tf(t,Y_t,Z_t)dt-\int_0^Tg(t,Y_t,Z_t)dB_t+\int_0^TZ_tdW_t
\end{eqnarray*}
and the estimates (5) and (6), we show the following inequalities,
\begin{eqnarray*}
E[K_T^2]&\leq&
cE[\xi^2+\int_0^Tf^2(s,0,0)^2ds+\int_0^Tg^2(s,0,0)ds+2\int_0^TS_sdK_s]\\
&\leq
&cE[\xi^2+\int_0^Tf^2(s,0,0)ds+\int_0^Tg^2(s,0,0)ds]+2c^2E(\sup_{0\leq
t\leq T}(S_t^+)^2)+\frac12E(K_T^2)\\
&\leq &
cE[\xi^2+\int_0^Tf^2(s,0,0)ds+\int_0^Tg^2(s,0,0)ds]+E(\sup_{0\leq
s\leq T}(S_s^+)^2)
\end{eqnarray*}
It follows easily that for each $t\in [0,T]$,
\begin{eqnarray*}
E(Y_t^2+\int_0^T|Z_t|^2dt+K_T^2)\leq
CE[\xi^2+\int_0^Tf^2(s,0,0)ds+\int_0^Tg^2(s,0,0)ds+(\sup_{0\leq
s\leq T}(S_s^+)^2)]
\end{eqnarray*}
The result then follows easily from Burkholder-Davis-Gundy
inequality.\quad\quad $\Box$\\

{\bf Proposition 3.2.} Let $(\xi,f,g,S)$ and $(\xi',f',g,S')$ be
two parameters satisfying the assumptions (i)-(iv). Suppose
$(Y,Z,K)$ is a solution of the reflected BDSDE $(\xi,f,g,S)$ and
$(Y',Z',K')$ is a solution of the reflected BDSDE
$(\xi',f',g,S')$. Define
$$\triangle \xi=\xi'-\xi,\quad \triangle f=f-f',\quad \triangle S=S-S',\quad \triangle Y=Y-Y',\quad \triangle Z=Z-Z',\quad \triangle K=K-K'$$
Then there exists a constant $C$ such that,
\begin{eqnarray}\nonumber
& & E(\sup_{0\leq t\leq T}|\triangle Y_t|^2+\int_0^T|\triangle
Z_t|^2dt+|\triangle K_T|^2) \\
& \leq & CE[|\triangle\xi|^2+\int_0^T|\triangle f(t,Y_t,Z_t)|^2dt
+cE[\sup_{0\leq t\leq T}({\triangle {S_t^+})^2}]^{\frac
12}\Psi_T^{\frac12}
\end{eqnarray}
where
$$\Psi_T=E[\xi^2+\int_0^Tf^2(t,0,0)dt+\sup_{0\leq t\leq T}(S_T^+)^2+{\xi'}^2+\int_0^T{f'}^2(t,0,0)dt+\sup_{0\leq t\leq T}({S'}^+_T)^2]$$
{\bf Proof.} The computation are similar to those in the previous
proof, so we'll only sketch the argument. Since
$\int_t^T(\triangle Y_s-\triangle S_s)d(\triangle K_s)\leq 0$,
\begin{eqnarray*}
E|\triangle Y_t|^2+E\int_t^T|\triangle Z_s|^2ds&\leq &
E|\triangle\xi|^2+2E\int_t^T\triangle f(s,Y_s,Z_s)\triangle Y_sds\\
& & +2E\int_t^T(f(s,Y_s,Z_s)-f(s,Y'_s,Z'_s))\triangle
Y_sds\\
& &
+E\int_t^T(g(s,Y_s,Z_s)-f(g,Y'_s,Z'_s))^2ds+2E\int_t^T\triangle
S_sd(\triangle K_s)
\end{eqnarray*}
Arguments already used in the previous proof lead to
\begin{eqnarray*}
E|\triangle Y_t|^2+\frac{1-\alpha}{2}E\int_t^T|\triangle
Z_s|^2ds&\leq & c[E|\triangle\xi|^2+\int_t^T|\triangle
f(s,Y_s,Z_s)|^2ds+\int_t^T|\triangle
Y_s|^2ds\\
& & +(\sup_{0\leq t\leq T}|\triangle S_t|)(K_T+K'_T)]
\end{eqnarray*}
Using Gronwall's lemma, Proposition 3.1 and the
Burkholder-Davis-Gundy inequality, we obtain inequality
(7).\quad\quad $\Box$\\

From the Proposition 3.2, We deduce immediately the following
uniqueness result  when $\xi=\xi',\quad f'=f,\quad S=S'$.

{\bf Theorem 3.3.} Under the assumption (i)-(iv), there exists at
most one measurable triple $(Y,Z,K)$, which satisfies (v)-(viii).
\section{Existence of a solution of reflected BDSDE: approximation via penalization.}
In this section, we'll give the result of existence via
penalization which is slightly different from [KKPPQ].

For each $n\in$ {\bf N}, let $(Y^n,Z^n)$ denote the unique pair of
${\cal{F}}_t$ measurable processes with valued in $R\times R^d$,
satisfying
$$E\int _0^T|Z_t^n|^2dt<\infty$$
and
\begin{eqnarray}
Y_t^n=\xi+\int_t^Tf(s,Y_s^n,Z_s^n)ds+n\int_t^T(Y^n_s-S_s)^-ds+\int_t^Tg(s,Y_s^n,Z_s^n)dB_s-\int_t^TZ^n_sdW_s
\end{eqnarray}
where $\xi$, $f$ and $g$ satisfy the  assumptions stated in
Section 2. We define
$$K^n_t=n\int_0^t(Y^n_s-S_s)^-ds, \quad \quad 0\leq t\leq T$$
From [PP3], we get
$$E(\sup_{0\leq t\leq T}|Y^n_t|^2)< \infty .$$
We now establish a priori estimate, uniformly in $n$, on the
sequence $(Y^n,Z^n,K^n)$.
\begin{eqnarray}\nonumber
& & E|Y^n_t|^2+E\int_t^T|Z^n_s|^2ds\\
\nonumber &=&
E|\xi|^2+2E\int_t^Tf(s,Y^n_s,Z^n_s)Y^n_sds+E\int_t^Tg^2(s,Y^n_s,Z^n_s)ds+2E\int_t^TY^n_sdK^n_s\\
\nonumber
\end{eqnarray}
from (4), we obtain
\begin{eqnarray}\nonumber
E|Y^n_t|^2+E\int_t^T|Z^n_s|^2ds &\leq &
E|\xi|^2+2E\int_t^T(f(s,0,0)+K|Y^n_s|+|Z_s^n|)|Y^n_s|ds\\
\nonumber & &
+E\int_t^T[c(\alpha')(|Y^n_s|^2+g^2(s,0,0))+\alpha'|Z_s^n|^2]ds+2E\int_t^TS_sdK_s^n\\
\nonumber &\leq &
E|\xi|^2+E\int_0^Tf^2(s,0,0)ds+E\int_0^Tg^2(s,0,0)ds\\
\nonumber & & +(2K+c(\alpha')+K^2\beta)E\int_t^T|Y^n_s|^2ds +
(\alpha'+\frac{1}{\beta})E\int_t^T|Z_s^n|^2ds \\
\nonumber & &
+\beta E[\sup_{0\leq t\leq
T}(S_t^+)^2]+\frac{1}{\beta}E[(K_T^n-K_t^n)^2]
\end{eqnarray}
where $\beta$ is positive. But
\begin{eqnarray*}
K^n_T-K^n_t=Y^n_t-\xi-\int_t^Tf(s,Y^n_s,Z^n_s)ds-\int_t^Tg(s,Y^n_s,Z^n_s)dB_s+\int_t^TZ_s^ndW_s,
\end{eqnarray*}
hence
\begin{eqnarray*}
E[(K^n_T-K^n_t)^2]\leq
c[E(|Y^n_t|^2)+E|\xi|^2+1+E\int_t^T(|Y^n_s|^2+|Z^n_s|^2)ds]
\end{eqnarray*}
choosing $\beta$ enough large, such that
$\alpha'+\frac{1+c}{\beta}\leq \bar{\alpha}<1$, then
\begin{eqnarray*}
E(|Y^n_t|^2)+(1-\bar{\alpha})E\int_t^T|Z^n_s|^2ds\leq
c(1+E\int_t^T|Y^n_s|^2ds)
\end{eqnarray*}
it then follows from Gronwall's lemma that
\begin{eqnarray*}
E(|Y^n_t|^2)+E\int_t^T|Z^n_s|^2ds+E[(K_T^n)^2]\leq c,\quad \quad
n\in {\bf N}
\end{eqnarray*}
furthermore, from Burkholder-Davis-Gundy inequality, we deduce
that
\begin{eqnarray}
E(\sup_{0\leq t\leq T
}|Y^n_t|^2)+E\int_t^T|Z^n_s|^2ds+E[(K_T^n)^2]\leq c,\quad \quad
n\in {\bf N}
\end{eqnarray}
note that if we define
\begin{eqnarray*}
f_n(t,y,z)&=& f(t,y,z)+n(y-S_t)^-,\\
f_n(t,y,z)&\leq&  f_{n+1}(t,y,z),
\end{eqnarray*}
it follows from lemma 2.2 that $Y^n_t\leq Y^{n+1}_t$, $0\leq t\leq
T$, a.e. Hence
$$Y^n_t\uparrow Y_t,\quad \quad 0\leq t\leq T,\quad a.e.$$
and from (9) and Fatou's lemma,
$$E(\sup_{0\leq t\leq T}Y^2_t)\leq c.$$
It then follows by dominated convergence that
\begin{eqnarray}
E\int_0^T(Y_t-Y^n_t)^2dt\rightarrow 0,\quad as \quad n\rightarrow
\infty
\end{eqnarray}
Next, we'll prove $Z_t^n\rightarrow Z_n$ in ${\mathbf{M}}^2$.\\
Applying It$\hat{o}$'s formula to $(Y^n-Y^p)$ and the function
$y\rightarrow y^2$ ,
\begin{eqnarray*}
& & E(|Y^n_t-Y^p_t|^2)+E\int_t^T|Z^n_s-Z^p_s|^2ds \\&=&
2E\int_t^T[f(s,Y^n_s,Z^n_s)-f(s,Y^p_s,Z^p_s)](Y^n_s-Y^p_s)ds\\
& &
+E\int_t^T|g(s,Y^n_s,Z^n_s)-g(s,Y^p_s,Z^p_s)|^2ds+2E\int_t^T(Y^n_s-Y^p_s)d(K^n_s-K^p_s)\\
&\leq & 2KE\int_t^T(|Y^n_s-Y^p_s|^2+|Y^n_s-Y^p_s|\cdot
|Z^n_s-Z^p_s|)ds+KE\int_t^T|Y^n_s-Y^p_s|^2ds\\
 & &+\alpha E\int_t^T|Z^n_s-Z^p_s|^2ds
+2E\int_t^T(Y^n_s-S_s)^-dK^p_s+2E\int_t^T(Y^p_s-S_s)^-dK^n_s
\end{eqnarray*}
from $2ab\leq \frac{2}{1-\alpha}a^2+\frac{1-\alpha}{2}b^2$, then
\begin{eqnarray*}
& & E(|Y^n_t-Y^p_t|^2)+E\int_t^T|Z^n_s-Z^p_s|^2ds \\
&=&
(3K+K^2\frac{2}{1-\alpha})E\int_t^T|Y^n_s-Y^p_s|^2ds+\frac{1+\alpha}{2}E\int_t^T|Z^n_s-Z^p_s|^2ds\\
& & +2E\int_t^T(Y^n_s-S_s)^-dK_s^p+2E\int_t^T(Y^p_s-S_s)^-dK_s^n
\end{eqnarray*}
\begin{eqnarray*}
\frac{1-\alpha}{2}E\int_t^T|Z_s^n-Z_s^p|^2ds&\leq&
cE\int_t^T|Y^n_s-Y^p_s|^2ds+(E(\sup_{0\leq t\leq
T}|(Y^n_t-S_t)^-|^2)\cdot E(K^p_T)^2)^{\frac12}\\
& & +(E(\sup_{0\leq t\leq T}|(Y^p_t-S_t)^-|^2)\cdot
E(K^n_T)^2)^{\frac12}
\end{eqnarray*}
so
\begin{eqnarray}\nonumber
E\int_t^T|Z_s^n-Z_s^p|^2ds&\leq&
c[E\int_t^T|Y^n_s-Y^p_s|^2ds+(E(\sup_{0\leq t\leq
T}|(Y^n_t-S_t)^-|^2)\cdot E(K^p_T)^2))^{\frac12}\\
 & & +(E(\sup_{0\leq t\leq T}|(Y^p_t-S_t)^-|^2)\cdot
E(K^n_T)^2]))^{\frac12}
\end{eqnarray}
now, we give the proof that
\begin{eqnarray}
E(\sup_{0\leq t\leq T}|(Y^n_t-S_t)^-|^2)\rightarrow 0,\quad \quad
as \quad n\rightarrow \infty
\end{eqnarray}
Since $Y^n_t\geq Y^0_t$, we can w.l.o.g. replace $S_t$ by $S_t\vee
Y^0_t$;so assume that $E(\sup_{t\leq T}S^2_t)<\infty$. We first
want to compare a.s. $Y_t$ and $S_t$ for all $t\in [0,T]$, while
we do not know yet that $Y$ is a.s. continuous. From the
comparison theorem for BDSDE's, we have that a.s. $Y^n_t\geq
\tilde{Y}^n_t$, $0\leq t\leq T$, $n\in \mathbf{N}$, where
$\{\tilde{Y}^n_t,\tilde{Z}^n_t;0\leq t\leq T\}$ is the unique
solution of the BDSDE:
\begin{eqnarray*}
\tilde{Y}^n_t=\xi+\int_t^Tf(s,Y^n_s,Z^n_s)ds+n\int_t^T(S_t-\tilde{Y}^n_s)ds+\int_t^Tg(s,Y^n_s,Z^n_s)dB_s-\int_t^T\tilde{Z}_s^ndW_s
\end{eqnarray*}
Let $\nu$ be a stopping time such that $0\leq \nu\leq T$. Then
\begin{eqnarray*}
\tilde{Y}^n_t&=&E^{{\cal{F}}_{\nu}}[e^{-n(T-\nu)}\xi+\int_{\nu}^Te^{-n(s-\nu)}f(s,Y^n_s,Z^n_s)ds+n\int_{\nu}^Te^{-n(s-\nu)}S_sds]\\
             & &+\int_{\nu}^Te^{-n(s-\nu)}g(s,Y^n_s,Z^n_s)dB_s
\end{eqnarray*}
It is easily seen that
$$e^{-n(T-\nu)}\xi+n\int_{\nu}^Te^{-n(s-\nu)}S_sds\rightarrow \xi{\mathbf{1}}_{\nu=T}+S_{\nu}{\mathbf{1}}_{\nu<T}$$
a.s. and in ${\mathbf{L}}^2$, and the conditional expectation
converges also in ${\mathbf{L}}^2$. Moreover,
$$|\int_{\nu}^Te^{-n(s-\nu)}f(s,Y^n_s,Z^n_s)ds|\leq \frac{1}{\sqrt{2n}}(\int_0^Tf^2(s,Y^n_s,Z^n_s)ds)^{\frac12}$$
hence
$E^{{\cal{F}}_{\nu}}\int_{\nu}^Te^{-n(s-\nu)}f(s,Y^n_s,Z^n_s)ds\rightarrow
0$ in ${\mathbf{L}}^2$, as $n\rightarrow \infty$. \\
and
\begin{eqnarray*}
E(\int_{\nu}^Tg(s,Y^n_s,Z^n_s)dB_s)^2&\leq& cE\int_0^Te^{-2n(s-\nu)}g^2(s,Y^n_s,Z^n_s)ds\\
&\leq & \frac{c}{4n}E\int_0^Tg^4(s,Y^n_s,Z^n_s)ds\rightarrow 0
\end{eqnarray*}
Consequently, $\tilde{Y}^n_s\rightarrow
\xi{\mathbf{1}}_{\nu=T}+S_{\nu}{\mathbf{1}}_{\nu<T}$ in mean
square, and $Y_{\nu\geq S_{\nu}}$ a.s. From this and the section
theorem in Dellacherie and Meyer [DM], it follows that a.s.
$$Y^n_t\geq S_t,\quad 0\leq t\leq T$$
Hence $(Y^n_t-S_t)^-\searrow 0$, $0\leq t\leq T$, a.s., and from
Dini's theorem the convergence is uniform in $t$. The result
finally follows by dominated convergence, since $(Y^n_t-S_t)^-\leq
(S_t-Y^0_t)^+\leq |S_t|+|Y^0_t|$.\\
From above property, (11) and
(10), hence
$$E\int_0^T|Z^n_s-Z^p_s|^2ds\rightarrow 0,\quad E\int_0^T|Y^n_s-Y^p_s|^2ds\rightarrow 0,\quad as \quad n,p \rightarrow \infty$$
Now, we want to prove the process $Y$ is continuous. Similar to
above proof,
\begin{eqnarray*}
& & |Y^n_t-Y^p_t|^2+\int_t^T|Z^n_s-Z^p_s|^2ds \\&=&
2\int_t^T[f(s,Y^n_s,Z^n_s)-f(s,Y^p_s,Z^p_s)](Y^n_s-Y^p_s)ds+\int_t^T|g(s,Y^n_s,Z^n_s)-g(s,Y^p_s,Z^p_s)|^2ds\\
& &
+2\int_t^T[g(s,Y^n_s,Z^n_s)-g(s,Y^p_s,Z^p_s)](Y^n_s-Y^p_s)dB_s-2\int_t^T(Y^n_s-Y^p_s)(Z^n_s-Z^p_s)dW_s\\
& & +2\int_t^T(Y^n_s-Y^p_s)d(K^n_s-K^p_s)\\
\end{eqnarray*}
and
\begin{eqnarray*}
\sup_{0\leq t\leq T}|Y^n_t-Y^p_t|^2 &\leq &
2\int_0^T|f(s,Y^n,Z^n_s)-f(s,Y^p,Z^p_s)|\cdot|Y^n_s-Y^p_s|ds+2\int_0^T(Y^p_s-S_s)^-dK^n_s\\
& &
+\int_0^T|g(s,Y^n,Z^n_s)-g(s,Y^p,Z^p_s)|^2ds+2\int_0^T(Y^n_s-S_s)^-dK^p_s\\
& & +2\sup_{0\leq t\leq
T}|\int_t^Tg(s,Y^n,Z^n_s)-g(s,Y^p,Z^p_s)(Y^n_s-Y^p_s)dB_s|\\
& & +2\sup_{0\leq t\leq T}|\int_t^T(Y^n_s-Y^p_s)(Z^n_s-Z^p_s)dW_s|
\end{eqnarray*}
and from Burkholder-Davis-Gundy inequality and $2ab\leq \beta
a^2+\frac{1}{\beta}b^2$, we have
\begin{eqnarray*}
E\sup_{0\leq t\leq T}|Y^n_t-Y^p_t|^2 &\leq &
\frac12 E\sup_{0\leq t\leq T}|Y^n_t-Y^p_t|^2+cE\int_0^T(|Y^n_s-Y^p_s|^2+|Z^n_s-Z^p_s|^2)ds\\
& & +(E[\sup_{0\leq t\leq T}|(Y^n_t-S_t)^-|^2]\cdot
E|K_T^p|^2)^{\frac12}\\
& &+(E[\sup_{0\leq t\leq T}|(Y^p_t-S_t)^-|^2]\cdot
E|K_T^n|^2)^{\frac12}
\end{eqnarray*}
hence, $E(\sup_{0\leq t\leq T}|Y^n_t-Y^p_t|^2)\rightarrow 0$, as
$n,p\rightarrow \infty$.\\
from which we get $Y^n$ convergence uniformly in $t$ to $Y$, a.s.
and $Y$ is a continuous process.\\

Denote $K^n_t=n\int_0^t(Y^n_s-S_s)^-ds$, since $K^n_{\cdot}
\nearrow$ as $n\nearrow $, and from $E((K^n_T)^2)\leq C$, $\forall
n\in N$ we have $K^n_T\nearrow K_T$ and $E(K_T)^2\leq C$, that is
$K_T<
\infty$, a.s.\\
Since
\begin{eqnarray*}
|K_t^n-K_t^p|&\leq&
|Y_t^n-Y_t^p|+|Y_0^n-Y_0^p|+|\int_0^t(f(s,Y^n_s,Z^n_s)-f(s,Y^p_s,Z^p_s))ds|\\
& &
+|\int_0^t(g(s,Y^n_s,Z^n_s)-g(s,Y^p_s,Z^p_s))dB_s|+|\int_0^t(Z_s^n-Z_s^p)dW_s|
\end{eqnarray*}
\begin{eqnarray*}
E(\sup_{0\leq t \leq T}|K_t^n-K_t^p|^2)& \leq& c\{E\sup_{0\leq
t\leq T
}|Y^n_s-Y^p_s|^2+E|Y^n_0-Y^p_0|^2\\
& & +E\int_0^T(f(s,Y^n_s,Z^n_s)-f(s,Y^p_s,Z^p_s))^2ds\\
& & +E(\sup_{0\leq t\leq
T}|\int_0^tg(s,Y^n_s,Z^n_s)-g(s,Y^p_s,Z^p_s))dB_s|)\\
& & +E(\sup_{0\leq t\leq T}|\int_0^t(Z^n_s-Z^p_s)dW_s|\}
\end{eqnarray*}
We use the fact that $f$ and $g$ are Lipschitz functions, and the
Burkholder-Davis-Gundy inequality for the last terms, he obtain
$$E(\sup_{0\leq t\leq T}|K_t^n-K_t^p|^2)\rightarrow 0, \quad as \quad n,p \rightarrow \infty$$
consequently, there exists a pair $(Z,K)$ of measurable processes
which valued in $R^d\times R$, satisfying
$$E(\int_0^T(Z^n_t-Z^p_t)^2dt+\sup_{0\leq t\leq T}|K_t-K_t^n|^2)\rightarrow 0, \quad \quad as \quad n\rightarrow \infty$$
and (v), (vi) satisfied by the triple $(Y,Z,K)$ (obtained by
taking limit as $n\rightarrow \infty$), (vii) from (12). It
remains to
check that $\int_0^T(Y_t-S_t)dK_t=0$.\\
Clearly, $\{K_t\}$ is increasing. Moreover, we have just seen that
$(Y^n,K^n)$ tends to $(Y,K)$ uniformly in $t$ in probability. Then
the measure $dK^n$ tends to $dK$ weakly in probability,
$$\int_0^T(Y^n_t-S_t)dK_t^n\rightarrow \int_0^T(Y_t-S_t)dK_t,$$
in probability, as $n\rightarrow \infty$.\\
We deduce from the same argument and (12) that
$$\int_0^T(Y_t-S_t)dK_t\geq 0.$$
on the other hand, $$\int_0^T(Y^n_t-S_t)dK_t^n\leq 0,\quad \quad
n\in {\bf N}$$ hence,
$$\int_0^T(Y_t-S_t)dK_t= 0, \quad \quad a.s.$$
and we have proved that $(Y,Z,K)$ solves the reflected
BDSDE(2).\quad\quad $\Box$
\section{Comparison Theorem for reflectd BDSDE.}
We next give a comparison theorem, similar to that of [KKPPQ] and
[HLM] for reflected BSDEs.

{\bf Theorem 5.1.} Let $(\xi,f,g,S)$ and $(\xi',f',g,S')$ be two
sets of data, each one satisfying all the assumptions of (i)-(iv)
[with the exception that the Lipschtiz condition (H1) could be
satisfied by either $f$ or $f'$ only]. And suppose in addition the
following:

(1) $\xi\leq \xi'$, a.s.,

(2) $f(t,y,z)\leq f'(t,y,z)$, \quad $dP\bigotimes dt$, a.e.
$\forall (y,z)\in R\times R^d $,

(3) $S_t\leq S'_t$, $0\leq t\leq T$, a.s.\\
Let $(Y,Z,K)$ be a solution of the reflected BDSDE with data
$(\xi,f,g,S)$ and $(Y',Z',K')$ a solution of the reflected BDSDE
with data $(\xi',f',g,S')$. Then
$$Y_t\leq Y'_t,\quad 0\leq t\leq T, \quad a.s.$$
If $f$ and $f'$ all satisfy Lipschitz condition (iii), and $S=S'$,
then we also have $dK\geq dK'$, $P\_a.s.$

{\bf Proof.} Applying It$\hat{o}$'s formula to $|(Y_t-Y'_t)^+|^2$,
and taking expectation, we get
\begin{eqnarray*}
E|(Y_t-Y'_t)^+|^2+E\int_t^T{\mathbf{1}}_{\{Y_s>
Y'_s\}}|Z_s-Z'_s|^2ds\leq
2E\int_t^T(Y_s-Y'_s)^+[f(s,Y_s,Z_s)-f'(s,Y'_s,Z'_s))]ds\\
+E\int_t^T{\mathbf{1}}_{\{Y_s>
Y'_s\}}[g(s,Y_s,Z_s)-g(s,Y'_s,Z'_s)]^2ds+2E\int_t^T(Y_s-Y'_s)^+(dK_s-dK'_s)
\end{eqnarray*}
since on $\{Y_t>Y'_t\}$, $Y_t>S_t'>S_t$, then $dK_t=0$, so we have
$$\int_t^T(Y_s-Y'_s)^+(dK_s-dK'_s)=-\int_t^T(Y_s-Y'_s)^+dK'_s\leq 0$$
Assume now that the Lipschitz condition in the statement applied
to $f$, then
\begin{eqnarray*}
& & E|(Y_t-Y'_t)^+|^2+E\int_t^T{\mathbf{1}}_{\{Y_s>
Y'_s\}}|Z_s-Z'_s|^2ds\\
&\leq & 2E\int_t^T(Y_s-Y'_s)^+[f(s,Y_s,Z_s)-f'(s,Y'_s,Z'_s))]ds\\
& & +E\int_t^T{\mathbf{1}}_{\{Y_s>
Y'_s\}}[g(s,Y_s,Z_s)-g(s,Y'_s,Z'_s)]^2ds\\
&\leq &
2K\int_t^T(Y_s-Y'_s)^+[|Y_s-Y'_s|+|Z_s-Z'_s|]ds+E\int_t^T{\mathbf{1}}_{\{Y_s>
Y'_s\}}[K|Y_s-Y'_s|^2+|Z_s-Z'_s|^2]ds\\
&\leq &
(3K+K^2\frac{2}{1-\alpha})E\int_t^T|(Y_s-Y'_s)^+|^2ds+\frac{1+\alpha}{2}\int_t^T{\mathbf{1}}_{\{Y_s>
Y'_s\}}|Z_s-Z'_s|^2ds
\end{eqnarray*}
hence
$$E|(Y_t-Y'_t)^+|^2\leq \bar{K}E\int_t^T|(Y_s-Y'_s)^+|^2ds,$$
and from Gronwall's lemma, $(Y_t-Y'_t)^+=0$, $0\leq t\leq T$,
a.s.\\
If $f$ and $f'$ are all Lipschitz functions and $S=S'$, we
consider the following two BDSDEs:
\begin{eqnarray*}
Y^n_t&=&\xi+\int_t^Tf(s,Y_s^n,Z^n_s)ds+n\int_t^T(Y^n_s-S_s)^-ds+\int_t^Tg(s,Y^n_s,Z^n_s)dB_s-\int_t^TZ^n_sdW_s,\\
{Y'}^n_t&=&\xi'+\int_t^T{f'}(s,{Y'}_s^n,{Z'}^n_s)ds+n\int_t^T({Y'}^n_s-S_s)^-ds+\int_t^Tg(s,{Y'}^n_s,{Z'}^n_s)dB_s-\int_t^T{Z'}^n_sdW_s,
\end{eqnarray*}
from the comparison theorem of BDSDE [GS], we get $\forall n\geq
0$, $P\_a.s.$ $Y^n\leq {Y'}^n$. On the other hand, from the proof
of existence in section 4, we know that, $\forall t\in [0,T]$,
$P\_a.s.$,

(i) $Y^n_t\rightarrow Y_t$ (resp. ${Y'}^n_t\rightarrow Y'_t$), as
$n\rightarrow\infty$,

(ii) $K_t=\lim_{n\rightarrow \infty}n\int_0^t(Y_s^n-S_s)^-ds$, and
$K'_t=\lim_{n\rightarrow \infty}n\int_0^t({Y'}_s^n-S_s)^-ds$.\\
for $Y^n\leq {Y'}^n$, it follows that, $\forall s,r\in [0,T]$,
$K_s-K_r\geq K'_s-K'_r$ and $dK\geq dK'$, $P\_a.s$.\quad \quad
$\Box$
\section{Other results}
{\bf Lemma 6.1.} let $(Y,Z,K)$ be a solution of the above
reflected BDSDE, satisfying condition (vi) to (viii). Then for
each $t\in [0,T]$,
$$K_T-K_t=\sup_{t\leq u\leq T}(\xi+\int_u^Tf(s,Y_s,Z_s)ds+\int_u^Tg(s,Y_s,Z_s)dB_s-\int_u^TZ_sdW_s-S_u)^-$$
{\bf Proof.} The proof is similar to [KKPPQ]. Where
$((Y_{T-t}(\omega)-S_{T-t}(\omega)),
(K_T(\omega)-K_{T-t}(\omega)), 0\leq t\leq T)$ is the solution of
a Skorohod problem. Applying the Skorohod lemma with
$$x_t=(\xi+\int_{T-t}^Tf(s,Y_s,Z_s)ds+\int_{T-t}^Tg(s,Y_s,Z_s)dB_s-\int_{T-t}^TZ_sdW_s-S_{T-t})(\omega),$$
$k_t=(K_t-K_{T-t})(\omega)$,
$y_t=(Y_{T-t}-S_{T-t})(\omega)$.\quad\quad $\Box$

{\bf lemma 6.2.} Let $(Y,Z,K)'$ be a solution of the above
reflected BDSDE (2), satisfying (v)-(viii). Then for each $t\in
[0,T]$,
\begin{eqnarray}
Y_t=ess\sup_{\nu\in
\Gamma_t}\{E^{{\cal{F}}_t}[\int^{\nu}_tf(s,Y_s,Z_s)ds+S_{\nu}{\mathbf{1}}_{\{\nu<
T \}}+\xi{\mathbf{1}}_{\{\nu= T\}}+\int_t^{\nu}g(s,Y_s,Z_s)dB_s\}]
\end{eqnarray}
where $\Gamma$ is the set of all stopping times dominated by $T$,
and $\Gamma_t=\{\nu\in\Gamma;t\leq \nu\leq T\}$.

{\bf Proof.} Let $\nu\in\Gamma_t$,
\begin{eqnarray*}
Y_t&=&\xi+\int_t^{\nu}f(s,Y_s,Z_s)ds+K_{\nu}-K_t+\int_t^{\nu}g(s,Y_s,Z_s)dB_s-\int_t^{\nu}Z_sdW_s\\
   &=&E[\xi+\int_t^{\nu}f(s,Y_s,Z_s)ds+K_{\nu}-K_t|{\cal{F}}_t]+\int_t^{\nu}g(s,Y_s,Z_s)dB_s\\
   &\geq & E[\int_t^{\nu}f(s,Y_s,Z_s)ds+S_{\nu}{\mathbf{1}}_{\{\nu<T\}}
   +\xi{\mathbf{1}}_{\{\nu=T\}}|{\cal{F}}_t]+\int_t^{\nu}g(s,Y_s,Z_s)dB_s
\end{eqnarray*}
Now er choose an optimal element of $\Gamma_t$ in order to get the
reversed inequality. Let
$$D_t=\inf\{t\leq u\leq T;Y_u=S_u\}\wedge T$$
Now the condition $\int_0^T(Y_t-S_t)dK_t=0$ and the continuity of
$K$ imply that
$$K_{D_t}-K_t=0,$$
it follows that
$$Y_t=E[\int_t^{D_t}f(s,Y_s,Z_s)ds+S_{D_t}{\mathbf{1}}_{\{D_t<T\}}
   +\xi{\mathbf{1}}_{\{D_t=T\}}|{\cal{F}}_t]+\int_t^{D_t}g(s,Y_s,Z_s)dB_s$$
from above result, we get (13).\quad\quad $\Box$
\section{Two reflected BDSDE}
Assume

(i) $\xi\in \mathbf{L}^2$;\\
let the mappings $f:[0,T]\times\Omega\times R\times R^d$,
$g:[0,T]\times\Omega\times R\times R^d$ be jointly measurably and
such that

(ii) $\forall (y,z)\in R\times R^d$, $f(\cdot,y,z)\in
\mathbf{M}^2$, $g(\cdot,y,z)\in \mathbf{M}^4$;

(iii) there exist two constants $K>0$ and $0<\alpha<1$, $\forall
(y,z), (y'z')\in R\times R^d$,
$$
\left\{
\begin{array}{lcl}
|f(t,y,z)-f(t,y',z')|^2 &\leq & C(|y-y'|^2+|z-z'|^2)\\
|g(t,y,z)-g(t,y',z')|^2 &\leq &C|y-y'|^2+\alpha
|z-z'|^2\\
\end{array}
\right. \eqno{\mbox{(H1)}}
$$
and two obstacles $\{L_t\}$ and $\{U_t\}$, such that

(iv) $E(\sup_{t\leq T}(L_t^+)^2)\leq \infty$, $E(\sup_{t\leq
T}(U_t^-)^2)< \infty$, and $L_t\leq \xi\leq
 U_t$, $P\_a.s.$,
$L_t<U_t$ for all $0\leq t<T$, $P\_a.s.$\\
A solution of two reflected BDSDE is a measurable processes
$(Y,Z,K^+,K^-)$, valued in $R\times R^d\times R_+\times R_+$, such
that for $0\leq t \leq T$

(v) $Z\in \mathbf{M^2}$;

(v') $Y\in \mathbf{S}^2$, $K^+$, $K^-\in \mathbf{L}^2$;

(vi) $ Y_t=\xi+\int_t^Tf(s,Y_s,Z_s)ds+(K^+_T-K^+_t)-(K^-_T-K^-_t)
 +\int_t^Tg(s,Y_s,Z_s)dB_s-\int_t^TZ_sdW_s$

(vii) $L_t\leq Y_t\leq U_t$, $P\_a.s.$, for all $0\leq t\leq T$;

(viii) $\{K_t^+\}$, $\{K_t^-\}$ are continuous and increasing,
$K_0^+=K_0^-=0$, and
$$\int_0^T(Y_t-L_t)dK^+_t=\int_0^T(U_t-Y_t)dK^-_t=0,\quad P\_a.s.$$
We also need the following additional assumption (H2):\\
there exists a process
$$X_t=X_0-\int_0^tJ_sdW-s-V_t^++V_t^-,\quad \quad X_T=\xi$$
with $J\in\mathbf{M}^2$, $V^+$, $V^-$ are continuous and
increasing, s.t.
$$L_t\leq X_t\leq U_t,\quad P\_a.s.  \forall t\in[0, T], \quad L<U,\quad P\_a.s. \forall t\in [0,T) $$
\vskip 0.5cm We now divide several steps to prove the existence
and uniqueness of two reflected BDSDE.

Consider the following BDSDE, for any $n$, $m\geq 1$,
\begin{eqnarray}\nonumber
Y^{n,m}_t&=&\xi+\int_t^Tf(s, Y^{n,m}_s,Z^{n,m}_s)ds+m\int_t^T(L_s-
Y^{n,m}_s)^+ds-n\int_t^T( Y^{n,m}_s-U_s)^+ds\\
 & & +\int_t^Tg(s,
Y^{n,m}_s,Z^{n,m}_s)dB_s-\int_t^TZ^{n,m}_sdW_s
\end{eqnarray}
since $f(s,y,z)+m(L_t-y)^+-n(y-U_t)^+$ is Lipschitz in $(y,z)$
uniformly in $(t,\omega)$, Eq(14) has a unique solution, denoted
$(Y^{n,m},Z^{n,m})$. Then we have the follow priori estimates.

{\bf lemma 7.1.} There exists a constant $C$ independent of $n$,
$m$, s.t.
\begin{eqnarray*}
\sup_{t\leq
T}E(Y^{n,m}_t)^2+E(\int_0^T|Z^{n,m}_s|^2ds)+m^2E(\int_0^T(L_s-Y^{n,m}_s)^+ds)^2
+n^2E(\int_0^T(Y^{n,m}_s-U_s)^+ds)^2\leq C
\end{eqnarray*}

{\bf proof.} (1) Applying It$\hat{o}$'s formula to $Y^{n,m}$ and
$y\rightarrow y^2$, we get
\begin{eqnarray*}
E(Y^{n,m}_t)^2+E(\int_t^T|Z^{n,m}_s|^2ds)&=&
E(\xi)^2+2E\int_t^TY^{n,m}_sf(s,Y^{n,m}_s,Z^{n,m}_s)ds\\
& &+E\int_t^Tg^2(s,Y^{n,m}_s,Z^{n,m}_s)ds
+2mE\int_t^TY^{n,m}_s(L_s-Y^{n,m}_s)^+ds\\
& & -2nE\int_t^TY^{n,m}_s(Y^{n,m}_s-U_s)^+ds
\end{eqnarray*}
from $2ab\leq \beta a^2+\frac{1}{\beta}b^2$, we have
\begin{eqnarray}\nonumber
E(Y^{n,m}_t)^2+E(\int_t^T|Z^{n,m}_s|^2ds)&=&
E(\xi)^2+c(\alpha')E\int_t^Tg^2(s,0,0)ds+E\int_t^Tf^2(s,0,0)ds\\\nonumber
&
&+(2K+c(\alpha')+K^2\frac{2}{1-\alpha'})E\int_t^T|Y^{n,m}_s|^2ds+\frac{1+\alpha'}{2}E\int_t^T|Z^{n,m}_s|^2ds\\\nonumber
& & +\beta E(\sup_{s\leq
T}(L_s^+)^2)+\frac{1}{\beta}m^2E(\int_t^T(L_s-Y^{n,m}_s)^+ds)^2\\
& & +\beta E(\sup_{s\leq
T}(U_s^-)^2)+\frac{1}{\beta}n^2E(\int_t^T(Y^{n,m}_s-U_s)^+ds)^2
\end{eqnarray}
we use the fact that
$$(Y^{n,m}_s-L_s)(L_s-Y^{n,m}_s)^+\leq 0,\quad (Y^{n,m}_s-U_s)(Y^{n,m}_s-U_s)^+\geq 0$$
(2) We now prove there exists a constant $\bar{c}$ independent of
$n$, $m$, such that, for all $0\leq t\leq T$,
\begin{eqnarray}\nonumber
& &m^2E(\int_0^T(L_s-Y^{n,m}_s)^+ds)^2
+n^2E(\int_0^T(Y^{n,m}_s-U_s)^+ds)^2 \\&\leq &
\bar{c}(1+E\int_t^T|Y^{n,m}_s|^2ds+E\int_t^T|Z^{n,m}_s|^2ds)
\end{eqnarray}
for $0\leq t\leq T$, define
\begin{eqnarray*}
T_1=\inf(t\leq r\leq T, Y^{n,m}_r=U_r)\wedge T;\\
S_1=\inf(T_1< r\leq T, Y^{n,m_r=L_r})\wedge T;\\
T_2=\inf(S_1\leq r\leq T, Y^{n,m}_r=U_r)\wedge T;\quad \cdots
\mbox{and
so on.}\\
\end{eqnarray*}
Then $T_k\nearrow T$, $S_k\nearrow T$ as $k\rightarrow \infty$.\\
Since $L<U$ on $[0,T)$, we have $Y^{n,m}\geq L$ between $T_k$ and
$S_k$, so
\begin{eqnarray*}
Y^{n,m}_{T_k}&=&Y^{n,m}_{S_k}+\int_{T_k}^{S_k}f_n(s,Y^{n,m}_s,Z^{n,m}_s)ds-n\int_{T_k}^{S_k}(Y^{n,m}_s-U_s)^+ds\\
& &
+\int_{S_k}^{T_k}g(s,Y^{n,m}_s,Z^{n,m}_s)dB_s-\int_{T_k}^{S_k}Z_s^{n,m}dW_s
\end{eqnarray*}
on the other hand,
\begin{eqnarray*}
U_{T_k}=Y^{n,m}_{T_k}\geq X_{T_k} \quad\mbox{if}\quad T_k<T,\quad
Y^{n,m}_{T_k}=X_{T_k}=\xi \quad\mbox{if}\quad T_k=T
\end{eqnarray*}
\begin{eqnarray*}
Y^{n,m}_{S_k}=L_{S_k}\leq X_{S_k}\quad\mbox{if}\quad S_k<T, \quad
Y^{n,m}_{S_k}=X_{S_k}=\xi\quad\mbox{if}\quad S_k=T
\end{eqnarray*}
from above property, we get for all $k$,
\begin{eqnarray*}
n\int_{T_k}^{S_k}(Y^{n,m}_s-U_s)^+ds &\leq &
X_{S_k}-X_{T_k}+\int_{T_k}^{S_k}f(s,Y_s^{n,m},Z_s^{n,m})ds\\
& &
+\int_{T_k}^{S_k}f(s,Y_s^{n,m},Z_s^{n,m})dB_s-\int_{T_k}^{S_k}Z_s^{n,m}dW_s\\
&\leq &
\int_{T_k}^{S_k}|f(s,Y_s^{n,m},Z_s^{n,m})|ds+V^+_{S_k}-V^+_{T_k}+V^-_{S_k}-V^-_{T_k}\\
& &
+\int_{S_k}^{T_k}g(s,Y_s^{n,m},Z_s^{n,m})dB_s-\int_{T_k}^{S_k}(J_s+Z_s^{n,m})dW_s
\end{eqnarray*}
since between $S_k$ and $T_{k+1}$, $Y^{n,m}_s\leq U_s$, summing up
in $k$, we obtain
\begin{eqnarray*}
n\int_T^t(Y^{n,m}_s-U_s)^+ds &\leq &
\int_t^T|f(s,Y^{n,m}_s,Z^{n,m}_s)|ds+V^+_T-V_t^++V_T^--V_t^-\\
& &
+\int_t^Tg(s,Y^{n,m}_s,Z^{n,m}_s)(\sum_k{\mathbf{1}}_{[T_k,S_k)}(s))dB_s\\
& &
-\int_t^T(J_s+Z_s^{n,m})(\sum_k{\mathbf{1}}_{[T_k,S_k)}(s))dW_s
\end{eqnarray*}
Taking square and expectation, we get
\begin{eqnarray}\nonumber
& & n^2E(\int_T^t(Y^{n,m}_s-U_s)^+ds)^2+E(\int_t^T(J_s+Z_s^{n,m})(\sum_k{\mathbf{1}}_{[T_k,S_k)}(s))ds)^2\\
&\leq
&\frac12\bar{c}(1+E\int_t^T|Z_s^{n,m}|^2ds+E\int_t^T|Y_s^{n,m}|^2ds)
\end{eqnarray}
In the same way, we obtain
\begin{eqnarray}
m^2E(\int_T^t(L_s-Y^{n,m}_s)^+ds)^2\leq\frac12\bar{c}(1+E\int_t^T|Z_s^{n,m}|^2ds+E\int_t^T|Y_s^{n,m}|^2ds)
\end{eqnarray}
from (15) and (16), choosing $\beta=\frac{2\bar{c}}{\alpha'}$, we
get
$$E(Y^{n,m}_t)^2+\frac{1-\alpha'}{2}E(\int_t^T|Z^{n,m}_s|^2ds)\leq \bar{K}(1+E\int_t^T(Y^{n,m}_s)^2ds)$$
then from Gronwall's lemma,
$$\sup_{t\leq T}E(Y^{n,m}_t)^2\leq c,\quad E\int_t^T|Z_s^{n,m}|^2\leq c$$
so
$$m^2E(\int^T_t(L_s-Y^{n,m}_s)^+ds)^2\leq c,\quad n^2E(\int^T_t(Y^{n,m}_s-U_s)^+ds)^2\leq c\quad \quad \Box$$

we now introduce the follow one reflected BDSDE $(\xi,f,g,L)$ for
all integer.
\begin{eqnarray}\nonumber
Y^n_t &=& \xi+\int_t^Tf(s,Y^n_s,Z^n_s)ds+K_T^{n,+}-K^{n,+}_t-n\int_t^T(Y^n_s-U_s)^+ds\\
& & +\int_t^Tg(s,Y_s^n,Z_s^n)dB_s-\int_t^TZ^n_sdW_s
\end{eqnarray}
from Burkholder-Davis-Gundy inequality,
$$E(\sup_{0\leq t\leq T}(Y^{n,m}_t)^2)\leq c$$
As $m\rightarrow\infty$, $Y^{n,m}\nearrow Y^n$,
$m\int_0^T(L_s-Y^{n,m}_s)^+ds\nearrow K_T^{n,+}$,
$Z^{n,m}\rightarrow Z^n$ in $M^2$.[KKPPQ]\\
 where $(Y^n,Z^n,K^{n,+})$ is
the unique solution of Eq$(\xi,f,g,L)$. Then

{\bf Lemma 7.2.}
\begin{eqnarray}
E(\sup_{t\leq
T}(Y^n_t)^2)+E\int_0^T|Z^n_s|^2ds+E(K_T^{n,+})^2+n^2E\int_0^T(Y^n_s-U_s)^+ds)^2\leq
c
\end{eqnarray}
where the constant $c$ is independent of $n$.\\
For Eq$(\xi,f,g,L)$, we know $Y^n\geq L$, a.s. and from the
comparison theorem [GS] that $Y^n\searrow$, we conclude that there
exists a process $Y$ such that $Y^n\searrow Y$, and from Fatou's
Lemma,
$$E(\sup_{t\leq T}Y^2_t)\leq c$$
then follows by the dominated convergence theotem that,
$$E(\int_0^T(Y_t-Y^n_t)^2dt)\longrightarrow 0,\quad \mbox{as}\quad n\rightarrow \infty$$
Now if $n\geq p$, $Y^n\leq Y^p$, $dK^{n,+}\geq dK^{p,+}$ by the
comparison theorem [GS]. Then we want to prove that
$Z^n\rightarrow Z$ in $M^2$, as $n\rightarrow\infty$.

{\bf Lemmma 7.3.} $Z^n\rightarrow Z$ in $M^2$ as $n\rightarrow
\infty$.

{\bf proof.} we define
$\widetilde{K}_t^T=n\int_0^t(Y^n_s-U_s)^+ds$, $n>p$.\\
Applying It$\hat{o}$'s formula to $Y^n-Y^p$ and $y\rightarrow
y^2$,
\begin{eqnarray*}
E(Y^n_t-Y^p_t)^2+E\int_t^T|Z^n_s-Z^p_s|^2ds &=& 2E\int_t^T(Y^n_s-Y^p_s)(f(s,Y^n_s,Z^n_s)-f(s,Y^p_s,Z^p_s))ds\\
& &
+2E\int_t^T(Y^n_s-Y^p_s)(dK_s^{n,+}-dK^{p,+}_s)\\
& & -2E\int_t^T(Y^n_s-Y^p_s)(d\widetilde{K}_s^{n,+}-d\widetilde{K}^{p,+}_s)\\
& &+E\int_t^T(g(s,Y^n_s,Z^n_s)-g(s,Y^p_s,Z^p_s))^2ds \\
&\leq
&2E\int_t^T(Y^n_s-Y^p_s)(f(s,Y^n_s,Z^n_s)-f(s,Y^p_s,Z^p_s))ds\\
& &  -2E\int_t^T(Y^n_s-Y^p_s)(d\widetilde{K}_s^{n,+}-d\widetilde{K}^{p,+}_s)\\
& & +E\int_t^T(g(s,Y^n_s,Z^n_s)-g(s,Y^p_s,Z^p_s))^2ds \\
\end{eqnarray*}
where we have use
$$(Y^n_s-Y^p_s)(d\widetilde{K}_s^{n,+}-d\widetilde{K}^{p,+}_s)=[(Y^n_s-U_s)+(U_s-Y^p_s)](d\widetilde{K}_s^{n,+}-d\widetilde{K}^{p,+}_s)$$
and $(Y^n_s-U_s)d\widetilde{K}_s^{n,+}\geq 0$,
$(U_s-Y^p_s)d\widetilde{K}^{p,+}_s\leq 0$, we obtain
\begin{eqnarray*}
E(Y^n_t-Y^p_t)^2+\frac{1-\alpha}{2}E\int_t^T|Z^n_s-Z^p_s|^2ds &=&
cE\int_t^T(Y^n_s-Y^p_s)^2ds\\
& & +2E(\sup_{t\leq T}((Y^n_s-U_s)^+)^2)\cdot
E(p\int_t^T(Y^p_s-U_s)^+ds)^2\\
& &+2E(\sup_{t\leq T}((Y^p_s-U_s)^+)^2)\cdot
E(n\int_t^T(Y^n_s-U_s)^+ds)^2
\end{eqnarray*}
Now, we prove $E(\sup_{t\leq T}((Y^n_s-U_s)^+)^2)\rightarrow 0$,
as $n\rightarrow \infty$. \\
We consider the following Eq$(\xi,f(s,Y^n_s,Z^n_s)-n(y-U_s),g,L)$,
\begin{eqnarray*}
\hat{Y}^n_t &=& \xi+\int_t^Tf(s,Y^n_s,Z^n_s)ds+\widetilde{K}_T^{n,+}-\widetilde{K}^{n,+}_t-n\int_t^T(\hat{Y}^n_s-U_s)ds\\
& & +\int_t^Tg(s,Y_s^n,Z_s^n)dB_s-\int_t^T\hat{Z}^n_sdW_s
\end{eqnarray*}
which has a unique solution, denoted
$(\hat{Y}^n,\hat{Z}^n,\widetilde{K}^n)$. From the comparison
theorem [GS], $Y^n_t\leq \hat{Y}^n_t$, a.s. for all $0\leq t\leq
T$.\\
Then we have the following reflected BDSDE:
\begin{eqnarray*}
e^{-nt}\hat{Y}^n_t &=& e^{-nT}\xi+\int_t^Te^{-ns}f(s,Y^n_s,Z^n_s)ds+\int_t^Te^{-ns}d\widetilde{K}_t^{n,+}+n\int_t^Te^{-ns}U_sds\\
& &
+\int_t^Te^{-ns}g(s,Y_s^n,Z_s^n)dB_s-\int_t^Te^{-ns}\hat{Z}^n_sdW_s
\end{eqnarray*}
This process $\{e^{-nt}Y^n_t\}$ is the solution of the BDSDE with
the obstacle $\{e^{-nt}L_t\}$ with terminal value $e^{-nt}\xi$ and
coefficients $e^{-nt}f(t,Y^n_t,Z^n_t)+ne^{-nt}U_t$ and
$g(t,Y^n_t,Z^n_t)$.\\
Let $\nu$ be a stopping time such that $0\leq \nu \leq T$, then
\begin{eqnarray*}
\hat{Y}^n_{\nu}&=&ess\sup_{\tau\geq \nu}\{E(\xi
e^{-n(\tau-\nu)}{\mathbf{1}}_{\{\tau
=T\}}+L_{\tau}e^{-n(\tau-\nu)}{\mathbf{1}}_{\{\tau
<T\}}+n\int_{\nu}^{\tau}e^{-n(s-\nu)}U_sds)\\
& &
+\int_{\nu}^{\tau}e^{-n(s-\nu)}f(s,Y^n_s,z^n_s)ds|{\cal{F}}_{\nu})+\int_{\nu}^{\tau}e^{-n(s-\nu)}f(s,Y^n_s,z^n_s)dB_s\}\\
&\leq &
E(n\int_{\nu}^{\tau}e^{-n(s-\nu)}(U_s-X_s)ds|{\cal{F}}_{\nu})+E(\int_{\nu}^{\tau}e^{-n(s-\nu)}|f(s,Y^n_s,z^n_s)|ds|{\cal{F}}_{\nu})\\
& & +ess\sup_{\tau\geq\nu}E(n\int_{\nu}^{\tau}X_se^{-n(s-\nu)}ds+
e^{-n(\tau-\nu)}X_{\tau}{\mathbf{1}}_{\{\tau
<T\}}+e^{-n(\tau-\nu)}\xi{\mathbf{1}}_{\{\tau=T\}}|{\cal{F}}_{\nu})\\
&
&+ess\sup_{\tau\geq\nu}\int_{\nu}^{\tau}e^{-n(s-\nu)}g(s,Y^n_s,z^n_s)dB_s
\end{eqnarray*}
where we have use $L_t\leq X_t\leq U_t$, a.s. $t\in [0,T]$.\\
It is easily seen that
$$n\int_{\nu}^{\tau}e^{-n(s-\nu)}(U_s-X_s)ds\rightarrow (U_t-X_t){\mathbf{1}}_{\{\tau
<T\}},\quad n\rightarrow \infty$$ a.s. and in $\mathbf{L}^2$, and
the conditional expectation convergence also in $\mathbf{L}^2$.
and
$$\int_{\nu}^{\tau}e^{-n(s-\nu)}|f(s,Y^n_s,z^n_s)|ds\leq \frac{1}{\sqrt{2n}}(\int_0^Tf^2(s,Y^n_s,z^n_s)ds)^{\frac12}$$
hence,
$$E(\int_{\nu}^{\tau}e^{-n(s-\nu)}|f(s,Y^n_s,z^n_s)|ds|{\cal{F}}_{\nu})\rightarrow 0 $$
in $\mathbf{L}^2$ as $n\rightarrow\infty$. \\
moreover,
\begin{eqnarray*}
E(ess\sup_{\tau\geq\nu}(\int_{\nu}^{\tau}e^{-n(s-\nu)}g(s,Y^n_s,z^n_s)dB_s)^2)&\leq&
cE\int_0^Te^{-2n(s-\nu)}g^2(s,Y^n_s,z^n_s)ds\\
&\leq & cE(e^{-4n(s-\nu)}ds\cdot
\int_0^Tg^4(s,Y^n_s,z^n_s))^{\frac12}\\
&\leq &
\frac{c}{4n}E(\int_0^Tg^4(s,Y^n_s,z^n_s))^{\frac12}\rightarrow 0
\end{eqnarray*}
in $\mathbf{L}^2$, as $n\rightarrow \infty$. \\
Now consider the second term at the right of the above inequality,
since
$$e^{-n(\tau-\nu)}X_{\tau}+n\int_{\nu}^{\tau}e^{-n(s-\nu)}X_sds=X_{\nu}\int_{\nu}^{\tau}e^{-n(s-\nu)}dX_s$$
we have
\begin{eqnarray*}
& & ess\sup_{\tau\geq\nu}E(n\int_{\nu}^{\tau}X_se^{-n(s-\nu)}ds+
e^{-n(\tau-\nu)}X_{\tau}{\mathbf{1}}_{\{\tau
<T\}}+e^{-n(\tau-\nu)}\xi{\mathbf{1}}_{\{\tau=T\}}|{\cal{F}}_{\nu})\\
&=&ess\sup_{\tau\geq\nu}E(X_{\nu}+\int_{\nu}^{\tau}e^{-n(s-\nu)}dX_s|{\cal{F}}_{\nu})\\
&\leq &X_{\nu}{\mathbf{1}}_{\{\nu <T\}}+\xi {\mathbf{1}}_{\{\nu
=T\}}+E(\int_{\nu}^{\tau}e^{-n(s-\nu)}d(V^++V^-)_s|{\cal{F}}_{\nu})
\end{eqnarray*}
since
$E(\int_{\nu}^{\tau}e^{-n(s-\nu)}d(V^++V^-)_s|{\cal{F}}_{\nu})\rightarrow
0 $ in $\mathbf{L}^2$, as $n\rightarrow \infty$. We obtain finally
$$Y_{\nu}\leq \hat{Y}_{\nu}\leq U_{\nu}{\mathbf{1}}_{\{\nu
<T\}}+\xi {\mathbf{1}}_{\{\nu =T\}}\leq U_{\nu},\quad a.s.$$ From
above and the section theorem of Dellacherie and Meyer [DM], it
follows that, $Y_t\leq U_t$, $0\leq t\leq T$, a.s. \\
Hence $(Y^n_t-U_t)^+\searrow 0$, $0\leq t\leq T$, a.s., and from
Dini's theorem the convergence in uniformly in $t$. Then, the
result finally follows by the dominated convergence theorem, since
$(Y^n_t-U_t)^+\leq (Y^0_t-U_t)^+\leq |Y^0_t|+|U_t|$.\quad\quad
$\Box$

{\bf Lemma 7.4.} The process $Y$ is a continuous process.

{\bf Proof.} Let $n>p$. Using It$\hat{o}$'s formula to $Y^n-Y^p$
and $y\rightarrow y^2$,
\begin{eqnarray*}
(Y^n_t-Y^p_t)^2+\int_t^T|Z^n_s-Z^p_s|^2ds &=& 2\int_t^T(Y^n_s-Y^p_s)(f(s,Y^n_s,Z^n_s)-f(s,Y^p_s,Z^p_s))ds\\
& & -2\int_t^T(Y^n_s-Y^p_s)(d\widetilde{K}_s^{n,+}-d\widetilde{K}^{p,+}_s)\\
& &+\int_t^T(g(s,Y^n_s,Z^n_s)-g(s,Y^p_s,Z^p_s))^2ds \\
& &
+\int_t^T(Y^n_s-Y^p_s)(g(s,Y^n_s,Z^n_s)-g(s,Y^p_s,Z^p_s))dB_s\\
& & -\int_t^T(Y^n_s-Y^p_s)(Z^n_s-Z^p_s)dW_s\\
\end{eqnarray*}
then
\begin{eqnarray*}
E(\sup_{t\leq T}(Y^n_t-Y^p_t)^2)+E\int_t^T|Z^n_s-Z^p_s|^2ds &\leq & E\int_t^T|Y^n_s-Y^p_s|\cdot|f(s,Y^n_s,Z^n_s)-f(s,Y^p_s,Z^p_s)|ds\\
& &+\int_t^T|g(s,Y^n_s,Z^n_s)-g(s,Y^p_s,Z^p_s)|^2ds\\
& & +2E(\sup_{t\leq T}(Y^n_t-U_t)^+\cdot p\int_0^T(Y^p_s-U_s)^+ds\\
& & +2E(\sup_{t\leq T}(Y^p_t-U_t)^+\cdot n\int_0^T(Y^n_s-U_s)^+ds\\
& &+E(\sup_{t\leq T}\int_t^T|(Y^n_s-Y^p_s)(g(s,Y^n_s,Z^n_s)-g(s,Y^p_s,Z^p_s))dB_s|)\\
& &+E(\sup_{t\leq T}\int_t^T|(Y^n_s-Y^p_s)(Z^n_s-Z^p_s)dW_s|\\
\end{eqnarray*}
Using B-D-G inequality and $f$ and $g$ is uniformly Lipschitz in
$(y,z)$, we deduce,
$$E(\sup_{t\leq T}(Y^n_t-Y^p_t)^2)\rightarrow 0,\quad \mbox{as}\quad n,p\rightarrow \infty $$
from which we get that $Y^n$ convergence uniformly in $t$ to $y$,
$P\_a.s.$ and that $Y$ is continuous process.\quad\quad$\Box$

\end{document}